\documentclass{amsproc}
\usepackage{graphicx}
\usepackage{enumerate}

\makeatletter
\let\@cleartopmattertags\relax
\newcommand\articleend
 {\enddoc@text
  \let\authors\@empty
  \let\contribs\@empty
  \let\xcontribs\@empty
  \let\toccontribs\@empty
  \let\addresses\@empty
  \let\thankses\@empty
  \let\dedicatory\@empty
\newpage
 }
\let\@wraptoccontribs\wraptoccontribs
\makeatother

\newtheorem{theorem}{Theorem}[section]

\newtheorem{corollary}[theorem]{Corollary}

\theoremstyle{definition}

\theoremstyle{remark}
\newtheorem{remark}[theorem]{Remark}

\numberwithin{equation}{section}

\newcommand{\F}{\mathbb F}
\newcommand{\Z}{\mathbb Z}
\newcommand{\Q}{\mathbb Q}
\newcommand{\R}{\mathbb R}

\newcommand{\C}{\mathbb C}

\DeclareMathOperator{\Gal}{Gal}

\renewcommand{\phi}{\varphi}

\def\epsilon{\varepsilon}
\def\kappa{\varkappa}

\renewcommand{\pod}[1]{\if@display\mkern10mu\else\mkern6mu\fi(#1)}
\renewcommand{\pmod}[1]{\pod{{\operator@font mod}\mkern6mu#1}}

\usepackage[left=2.5cm,right=2.5cm, top=3cm,bottom=3cm,bindingoffset=0cm]{geometry}

\begin{document}
\allowdisplaybreaks

\title{Complete Description of Measures Corresponding to Abelian Varieties over Finite Fields}

\author{Nikolai S. Nadirashvili and Michael A. Tsfasman}

\address{CNRS and IITP}

\thanks{Supported in part by ANR project FLAIR (ANR-17-CE40-0012)}




\email{mtsfasman@yandex.ru; nnicolas@yandex.ru}

\date{2021}


\maketitle

\begin{abstract}
We study probability measures corresponding to families of abelian varieties over a finite field.
These measures play an important role in the Tsfasman--Vl\u{a}du\c{t} theory of
asymptotic zeta-functions defining completely the limit zeta-function of the family.
J.-P.~Serre, using results of R.M.~Robinson on conjugate algebraic integers, 
described the possible set of measures than can correspond to families of abelian
varieties over a finite field. The problem whether all such measures actually occur was left open. 
Moreover, Serre supposed that not all such measures correspond to abelian varieties (for example, the Lebesgue measure on a segment).
Here we settle Serre's problem proving that Serre conditions are sufficient, and  thus describe completely the set of measures 
corresponding to abelian varieties.
\end{abstract}

\section{Preliminaries}

\subsection{Abelian varieties over a finite field}

Let us recall the notation and background from \cite{Se} and \cite{Ts0}.
We fix once and for all a finite field $\F_q$.
Given an algebraic variety $A$ of dimension $g$ over it, consider
the multiset $\Omega=\{\omega_j\}$ of its inverse Frobenius roots.
These $2g$ roots are
\begin{itemize}
\item
pairwise complex-conjugate (real roots having even multiplicity);
\item
of absolute value $\sqrt q$ ;
\item
roots of a monic polynomial with integer coefficients (i.e., they are algebraic
integers and if $\omega\in \Omega$ then all its $\Gal (\bar\Q / \Q)$-conjugates are
also in $\Omega$).
\end{itemize}

Call the multiset $\{\omega_j\}$ with these properties an \emph{integer Weil system} or just a \emph{Weil
system}. The Weil system determines the zeta-function of $A$.

Instead of considering Weil systems $\Omega=\{\omega_i\}$ on the circle we can sum pairs of complex conjugate ones and consider 
$X=\{x_i\}, i=1, ... , g$, where $x_i= \omega_i + \bar \omega_i$.

To the Weil system $\Omega$ there corresponds the probability measure
$$
\mu_\Omega= {\frac{1}{g}}\sum_{x_i \in X}\delta_{x_i},
$$
where $\delta_{x_i}$ is the Dirac skyscraper measure at $x_i$. We work in
the weak topology space of measures dual to continuous functions on the segment $I=[-2\sqrt q, 2\sqrt q]$, writing the pairing
as $\langle f,\mu \rangle = \int f d\mu$.

We are interested in the corresponding asymptotic problem.

For a family of Weil systems $F= \{\Omega_j\} $ of growing cardinality $g_i$
we consider the limit in the weak topology
$$
\mu_F =\lim_j \mu_{\Omega_i}.
$$
Of course, it need not exist. When it exists, we say that the
family is \emph{asymptotically exact}. Such is for example any tower of embedded abelian varieties. Each
family contains an asymptotically exact subfamily, since the space of positive measures of given mass
on the circle being compact. This notion is compatible with that of asymptotic
exactness for curves and abelian varieties used in \cite{Ts,Ts/Vl1,Ts/Vl2}. 

We face the following natural questions, cf. \cite{Ts0}:

\begin{enumerate}[Q1.]

\item \emph{Describe the set of all limit measures corresponding to asymptotically exact
families of abelian varieties over\/ $\F_q$.}
\item \emph{Describe the set of all limit measures corresponding to asymptotically exact
families of curves over\/ $\F_q$.}

\end{enumerate}

The second question is treated in \cite{Ts}, \cite{Ts/Vl1}, \cite{Ts0}, and --- together with its analog for number fields --- in \cite{Ts/Vl2}. The answer {\it {a priori}} looks rather strange. In particular, all limit measures have continuous density. But this is specific for curves. 
Moreover, for curves --- not even to speak about number fields --- the complete answer is still unknown (see \cite{Ts0}).

For abelian varieties the situation is quite different. For example, the limit measure can be atomic.

\subsection{Algebraic integers}

By the Honda--Tate theorem \cite{Ta} for any $\Omega$ there exist an abelian variety over $\F_q$ whose set of inverse Frobenius roots is a multiple of $\Omega$. (The word multiple here is important, $\Omega$ itself is not always realisable.)

The set $\Omega=\{\omega_j\}$ of cardinality $2g$ uniquely determines and is uniquely determined by the set $\bar\Omega=\{\omega_j+\bar\omega_j\}$ of cardinality $g$,  where for each pair of conjugate elements of $\Omega$ we take their sum. $\bar\Omega$ is a set of $g$ totally real algebraic numbers on the segment $I=[-2\sqrt{q},2\sqrt{q}]$. Thus the description of possible sets 
$\Omega$ becomes the question in algebraic number theory:
\begin{itemize}
\item
describe algebraic intergers such that all their $\Gal (\bar\Q / \Q)$-conjugates are
 in  $[-2\sqrt{q},2\sqrt{q}]$.
\end{itemize}

For each $\Gal (\bar\Q / \Q)$-stable multiset $X=\{x_i\}, i=1, ... , g$ of totally real algebraic numbers lying on the given segment $I$ we define the probability measure (i.e., measure of total mass 1), normalizing the sum of atomic measures at these numbers, $\mu_X = {\frac{1}{g}}\sum_{i=1}^g \delta_{x_i}$ . Consider a sequence ${\bf X} = \{X_j\}$ of such sets with growing $g_j$ and the limit measure 
$\mu_{\bf X} = \lim {\frac {1}{g_j}} \mu_{X_j}$. Of course, such measure need not exist. If the limit exists, we say that the family $\bf X$ is \emph{asymptotically exact}. (For us a \emph{family} means just a \emph{sequence with growing} $g_j$.)

Any family contains an  asymptotically exact subfamily. Any tower (family with $X_j \supset X_{j-1}$) is asymptotically exact.

We come to the following problem in algebraic number theory, equivalent to the above question Q1.
 
\begin{enumerate}[Q3.]
\item
Describe the set of all limit probability measures corresponding to asymptotically exact
families of $\Gal (\bar\Q / \Q)$-stable multisets of totally real algebraic integers lying on the given segment $I$.

\end{enumerate}

If we take some $\Gal (\bar\Q / \Q)$-stable multiset  $X$ and set $X_n=n X, n\to \infty$, then the limit is just $\mu_X$ which is an atomic measure. (This corresponds to the family of powers of a given abelian variety.)

Any atomic $\Gal (\bar\Q / \Q)$-stable measure supported in totally real algebraic integers on $I$ can be obtained in this way.

Recall that every positive measure $\mu$ on a segment is a sum of $\mu_a$ and $\mu_{d}$, where $\mu_a$ is atomic and 
$\mu_{d}$ is atom-free (or diffuse), i.e., $\mu_{d}(x)=0$ for any point $x$.

So our question is reduced to the following

\begin{enumerate}[Q4.]
\item
Describe the set of all limit  \emph{diffuse} probability measures corresponding to asymptotically exact
families of $\Gal (\bar\Q / \Q)$-stable multisets of totally real algebraic integers lying on the given segment $I$.

\end{enumerate}

\subsection{Capacity and equilibrium measure}
Following \cite{Se} we introduce some notation.
Consider a compact $K\subset\C$. For a positive measure $\mu$ on $K$ define
$$
I(\mu)=\int\int_{K\times K} \ln\vert x-y \vert \mu (x) \mu (y).
$$
Note that $I(\mu)\in \R\cup\{-\infty\}$.
Set $v(K)=\sup_{\mu} I(\mu) \in \R\cup\{-\infty\}$, where $\mu$ runs over all positive probability measures  supported in $K$.
Then $v(K)$ is called the  \emph {log-capacity} of $K$, and the \emph{capacity} of $K$ is defined as
${\rm{cap}}(K)=e^{v(K)}$.

For a non-compact $Y\subset\C$ the log-capacity is defined as supremum of log-capacities of its compact subsets.

For a compact $K$ with log-capacity $v(K)>-\infty$ there exists a unique positive probability measure  $\mu$, such that $I(\mu) =v(K)$. This measure is diffuse. Its support can be smaller than $K$, but for $K\subset \R$ their difference has log-capacity $-\infty$. This measure $\mu=\mu_K$ is called the \emph{equivilibrium measure} of $K$.

\subsection{Algebraic integers with all conjugates in a given compact}
For a set $V\subset \C$  let ${\rm{Irr}}_V$ be the set of irreducible monic polynomials in $\Z[x]$ of degree at least $1$ such that all their roots lie in $V$. For such a polynomial $P(x)$ of degree $g$ let $\mu_P$ be the corresponding probability measure supported in its roots, $\mu_P = {\frac{1}{g}}\sum_{i=1}^g \delta_{x_i}$.
Now let $K\subset \C$ be compact. There are two quite different cases (Fekete, Szeg\"{e}, Robinson, cf. \cite{Se}), depending on the capacity of $K$.
\begin{itemize}
\item
If  $v(K) < 0$, then $\rm{Irr}_K$ is finite.
\item
If $K\subset \R$  is a union of finite number of segments and  $v(K)> 0$, then $\rm{Irr}_K$ is infinite.
\item
If $K\subset \C$ is complex conjugation stable and  $v(K)\ge 0$, then for any open $U\supset K$ the set $\rm{Irr}_U$ is infinite.
\end{itemize}

\subsection{Space of measures corresponding to integer polynomials}
Consider probability measures $\mu_P$, corresponding to irreducible polinomials with all roots in $K$. The set of such measures is either finite, or  enumerable; introduce any order on it and let ${\bf M}_n$ be the closure in the weak topology of the convex envelope of all such measures except for the first $n-1$ ones. Then 
$$
{\bf M}_1\supset  {\bf M}_2\supset{\bf M}_2\supset  \dots 
$$
Let ${\bf M}_\infty = \bigcap_n {\bf M}_n$. Suppose that $\rm{Irr}_K$ is infinite, then we know \cite{Se} : 
\begin{itemize}
\item
The set ${\bf M}_\infty$ is non-empty, convex, and compact.
\item
Any measure $\mu\in {\bf M}_\infty$ is diffuse.
\item
The log-capacity $v(\rm{Supp}(\mu))$ of its support is greater than or equal to $0$. 
\item
If $v(\rm{Supp}(\mu))=0$ then $\mu=\mu_K$ is the equilibrium measure of $K$.
\end{itemize}

\subsection{Potential}
Potential $ p_{\mu} : \C \to \R\cup \{-\infty\}$ of the measure $\mu$ is defined as
$$
p_{\mu}(z)= \int_{\C} \ln\vert w-z\vert \mu(w).
$$

Potential $p_K$ of the equilibrium measure $\mu_K$ of $K$ is everywhere in $\C$  greater than or equal to the logarithmic capacity $v(K)$. Moreover, it equals $v(K)$ on $K$ and is harmonic outside $K$.

\section{Criterion}
\subsection{Main result}

We are going to prove the following

\begin{theorem}{\label{main}}
Let $I$ be a segment $[a,b]\subset \R$, and let $\mu$ be a positive probability meausre  whose support lies in $I$.
Then the following conditions are equivalent:
\begin{enumerate}
\item
$\mu$ lies in the weak topology closure of equilibrium measures $\mu_K$ for compact $K\subseteq I$ of log-capacity $v(K) \ge 0$, i.e. ${\rm{cap}}(K)\ge 1$;
\item
$\mu$ lies in the weak topology closure of equilibrium measures $\mu_K$ for compact $K\subseteq I$ of log-capacity $v(K) = 0$;
\item
$p_{\mu}(z)\ge  0$ for all  $z\in \C$.
\end{enumerate}
\end{theorem}

\subsection{Application to ${\bf{M}}_\infty$}

Recall the theorem of Serre (\cite{Se}, thm 1.6.3), following from the above results of Robinson (see 1.4).

\begin{theorem}
Let $K\subset \R$ be a finite union of segments with $v(K)> 0$, let $E\subset K$ be a compact subset with $v(E)\ge 0$, and $\mu_E$ be its equilibrium measure. Then  $\mu_E \in {\bf{M}}_\infty$.

\end{theorem}

From this theorem and Theorem \ref{main} we immediately derive the following criterion.

\begin{theorem}
Let $K\subset \R$ be a segment with $v(K)>0$, then a positive probability  measure $\mu$  on $K$  lies in ${\bf{M}}_\infty$ if and only if its potential is non-negative, i.e., $p_{\mu}(z)\ge 0$ for all $z\in \C$. 

\end{theorem}

\begin{corollary}
For a prime power $q$, a positive difuse probability measure  on the segment $I_q=[-2\sqrt {q}, 2\sqrt {q}]$ is the limit measure for an asymptotically exact family of abelian varieties over $\F_q$ if and only if 
$p_{\mu}(z)\ge  0$ for all  $z\in \C$.
\end{corollary}

\begin{remark}
In particular, we get an affirmative answer to the question (1.6.8) of Serre \cite{Se}: the normalized Lebesgue measure of mass $1$ on a segment $I=[-a, a]$ lies in ${\bf M}_\infty={\bf M}_\infty(I)$ for any $a\ge e$.
Indeed, for the Lebesgue measure supported on the segment $[-a,a]$ the potential  is
$$
p_L (y) = (y+a)(\ln{(y+a)} - 1) + (a-y)(\ln{(a-y)} - 1).
$$
It reaches its minimum at $y=0$ and $p_L (0) = 2 a (\ln a - 1)$, which is positive for $a\ge e$.
In particular (since $2\sqrt 2 > e$), for any prime power $q$ the Lebesgue measure supported on $I_q$ is realizable for asymptotically exact families of abelian varieties over $\F_q$.
\end{remark}

\section{Generalities on potential}

We set ${\mathbb C }={\mathbb R }^2$, $z=x+i y$.  By $D_r\subset {\mathbb R }^2, \, r>0,$ we denote the disk $ |z|< r $, and by $D^x_r$  the disk $ |z - x|< r $; by  $S_a\subset {\mathbb R }^2, \, a>0$, the slab $-a<y<a$; by  $R(a,b,c)\subset {\mathbb R }^2, \, c>0,$ the rectangle $a<x<b, \, -c<y<c$; by
 $I\subset {\mathbb R }^2 $ the segment $[ (-Q,0),(Q,0)] $. 
 
 Let $E{\subset \mathbb R }^2$ be compact. Denote $u_E(z)=p_E(z) - v(E)$, where $p_E(z)$ is the logarithmic potential of the equilibrium measure of the set $E$. We call function $u_E(z)$ the normalised potential of the set $E$. Define harmonic functions $u_E^R(z), \, R>Q$ in $D_R\setminus E$ satisfying the following boundary conditions: $u_E^R(z)=0$ on the regular points of $\partial E$ and $u_E^R(z)=1$ on  $\partial D_R$.  Then $(\ln R )u_E^R(z)$ tends to $p_E(z)=u_E(z)+ v(E)$ as $R\to \infty$. 
 
  The basic properties of the log-capacity and of the solution of the exterior Dirichlet problem are well known 
 see, e.g. \cite{Do}, cf.  
 (Part 1, and especially XIII.18)
  
(1) Let $E\subset I$ and $(-a,a)\subset E$. Then $u_E^R(z) < C(\ln |z|- \ln a) $ for $|z|>2a$, where $C>0$ is an absolute constant.

Let $\mu$ be a probability measure on $I$, and  
$$p_{\mu}(z)= \int_{\C} \ln\vert w-z\vert \mu(w)$$ 
its potential.
The following properties of   $p_\mu$ are easy to check:

(2)  $p_\mu(z)$ is subharmonic. It is harmonic on ${\mathbb R}^2\setminus \rm{Supp}(\mu)$.

(3) $p_\mu(z)$ is an even function of $y$.

(4) $p_\mu(z)$ is a monotonically increasing function of $y$ for $y>0$ and a monotonically decreasing function of $y$ for $y<0$.

(5) $p_\mu(z) = \ln |z| +o(1)$ as $|z| \to \infty $,
and $\mu=2\pi \Delta p_{\mu}(z)$. 

(6) Assume that $p_\mu (z) \geq 0$. Since  $p_\mu (z) > 1$ on ${\mathbb R }^2\setminus D_r$, if $r>Q+e$, then (by the Harnack inequality for harmonic functions)  for any $a>0$ there exists $c(a)>0$ such that $p_\mu(z) > c(a)$ in ${\mathbb R }^2\setminus S_a$.

(7) Let $\{\mu_n\}$ be a sequence of probability measures on $I$. Then, $\mu_n\to \mu$ in the weak topology if and only if $p_{\mu_{n}}(z)\to  p_{\mu}(z)$ pointwise.

(8)  Let $E \subset {\mathbb R }^2$ be a compact nonempty set with a smooth boundary.  Let $f$ be a continous function on $\partial E$. Denote $M^- =\inf f, \, M^+= \sup f$. Choose a constant $m$ such that $M^- \leq m\leq M^+$. For any large enough $R$ we can solve the Dirichlet problem:
$$ \Delta u_R=0 \; {\rm in} \; D_R\setminus E,$$
$$u_R=f \; {\rm on} \; \partial E, \;\; u_R=m \; {\rm on} \; \partial D_R .$$ 
Choose a convergent sequence $u_{r_n}, \; r_n \to \infty$  in $ {\mathbb R }^2 \setminus E$. Then the function
$$u=\lim u_{r_n}$$
is a solution of the external Dirichlet problem
$$ \Delta u=0 \; in \;  {\mathbb R }^2\setminus E,$$
$$u=f \; {\rm on} \; \partial E .$$
By the maximum principle each of the functions $u_{r_n}$ satisfies the inequalities
$$M^-\leq u_{r_n}\leq M^+ ,$$
and hence the function $u$ satisfies the same inequalities
$$M^-\leq u\leq M^+ .$$

We show now that the last Dirichlet problem has a unique bounded solution. Assume by contradiction that for a given function $f$ there are two different bounded solutions $u_1$ and $u_2$. Then $u=u_1-u_2$ is a bounded harmonic function
vanishing on $\partial E$.

 Assume now that $D_{r_0} \subset E \subset D_{r_1}$. Let $v_R$ be a solution of the Dirichlet problem
$$ \Delta v_R=0 \; {\rm in} \; D_R\setminus E,$$
$$v_R=0 \; {\rm on} \; \partial E, \;\; v_R=\ln (R/r_0) \; {\rm on} \; \partial D_R$$ 
Then $\ln |x|/r_1 < v(x) < \ln |x|/r_0$ and we can choose a sequence $r_n\to \infty$ such that $v_{r_n} \to u$ in
$ {\mathbb R }^2 \setminus E$. Then $v$ is a positive harmonic function in  $ {\mathbb R }^2 \setminus E$,
vanishing on $\partial E$ and tending to infinity at infinity. Thus for any positive $\epsilon$ for sufficiently large $R$
in $D_R\setminus E$ by the maximum principle we have $u < \epsilon v$. Since the last inquality
holds for any positive $\epsilon $  we get $u=0$.

(9) As an immediate consequence of (8) we have the following maximum principle for  harmonic functions in the unbounded domains in dimension $2$.  Let $E \subset {\mathbb R }^2$ be compact, and let $u$ be a bounded harmonic function in ${\mathbb R }^2 \setminus E$
continuous in its closure. Then the following maximum principle for the exterior domain holds:
$$ \inf_{\partial E} u< u <\sup_{\partial E} u  . $$

(10) Let $s=[-\epsilon, \epsilon ]$ be a segment on the axis $x$, and $u$ be a harmonic function in $D_1\setminus s$
such that $u=0$ on $s$ and $u=1$ on $\partial D_1$. Then $u < C\ln K / \ln (1/\epsilon) $ in $D_{K\epsilon} $ for any $1<K<1/\epsilon$, here $C>0$ is an absolute constant.

(11) Let $x_1,\dots , x_n \in [0,1], \; l_1,\dots , l_n<\epsilon $.  Let $X$ be the union of segments
$$ X= \cup_{i=1}^n [x_i-l_i,x_i+l_i]$$
and
$$ U_K= \cup_{i=1}^n D_{Kl_i}^{x_i},$$
where $K>1$. Assume $v(X)<0$. Since $p_X(z) < \ln (1+R)$ on $D_R, R>2$ then there is an absolute constant $C>0$ such that 
$p_X(z)<v(X)- Cv(X)\ln K$ on $U_K$. Now, from (10) we get
$$ v(U_K) <v(X)- Cv(X)\ln K .$$

\section{Proof of Theorem\:\ref{main}}
Any measure $\mu$ can be appoximated by a sequence $\mu_n$ of atomic measures. Our goal is to prove that if for all $z$ the potential $p_{\mu}(z)\ge 0$, then for every large enough $n$ we can approximate these
$\mu_n$ by  $\mu'_n$, where $\mu'_n$ is the equilibrium measure of a set $S_n$ which is a union of  small segments containing ${\rm{Supp}}(\mu_n)$. 

(12) Let $\mu_n$ be a sequence of probability measures, each supported in a finite set of points of $I$, $I = [a, b]$, such that $\mu_n$
tends to $\mu$ weakly. Then $p_{\mu_n }(z)\to p_\mu(z)$ in ${\mathbb R }^2\setminus S_\delta$ uniformly for any $\delta >0$.

(13) Assume that $p_\mu(z) \geq 0$. Denote $G_n:=\{ z\in {\mathbb R}^2, p_{\mu_n}(z)<0\}$.

From (12) it follows that for any $\epsilon >0$
there exists some $N=N(\epsilon )$ such that  $G_n\subset S_\epsilon$ for any $n>N$.

(14) Let $i_a, 0\leq a\leq 0$ be the segment $[-a, a], y=0$, $v(z)$ be a harmonic in $D_1 \setminus i_a$ function, $v(z)=0$
on $i_a$ and $v(z)=1$ on $\partial D_1$. Then $v(z)< C\ln (|z|/a)/\ln (1/|a|)$ in $D_1 \setminus D_{2a}$ where $C>0$ is an absolute constant.

(15) The set $I_n=G_n\cap I$ is a finite union of segments, which we denote $[ (a_i^n,0),(b_i^n,0)] $. Let $p_{I_n}(z)$ be the 
potential of the equilibrium measure of the set $I_n$. In the following sections we prove that $u_{I_n}(z)$ tends pointwise to $p_\mu(z)$ 
as $n$ goes to infinity. Obviously it is sufficient to prove that
$$| p_{\mu_n}(z)-u_{I_n}(z)|\to 0 \; .$$
 Let $H_n\subset D_R$ be a sequence of open sets  such that $I_n \subset H_n$.
 Note that by (9) the result will follow if we prove that there is a sequence $\epsilon_n >0, \, \epsilon_n \to 0$ such that
$$| p_{\mu_n}(z)-u_{I_n}(z)|< \epsilon_n \;\; {\rm on}\; \; \partial H_n \; .$$

 Note that since $u_{I_n}(z)\geq 0$, we have $u_{I_n}(z) -p_{\mu_n}(z) \geq 0$ on $I_n$ and hence everywhere. If this is proved, then we have 
 $u_{I_n}(z) -p_{\mu_n}(z) \to 0$ on an open set, hence from the last inequality, from (9), and from the Harnack inequality
 $u_{I_n}(z) -p_{\mu_n}(z) \to 0$ pointwise everywhere.

We still have to choose the sets $H_n$ and to prove the above inequality on them.

 (16) Fix $k$. Let $G_k \subset S_{\epsilon}$. Then
$$ G_k \subset \cup_i R(a_i^k,b_i^k, \epsilon) \; .$$

Set $l_i=b_i^k-a_i^k$. 

For $l_i<\epsilon$ we want to prove the inequality 
$$p_{\mu_k}(z)> -C_1e^{-C_2A} \;\; {\rm on}  \;\; \partial R(a_i^k,b_i^k,Al_i),$$
where $C_1,C_2$ are absolute constants, $A>1$. 

Indeed, denote $G^i=G_k \cap R(a_i^k,b_i^k,\epsilon)$. Then $G^i$ is a compact subdomain of the
rectangle $R(a_i^k,b_i^k,\epsilon)$ and $p_{\mu_k}(z)=0$ on $\partial G^i$. Since the equilibrium measure $\mu_k$ of the segment
$[ (a_i^k,0),(b_i^k,0)] $ is less than or equal to $1$, it follows that  $$p_{\mu_k}(z)> -C_0 \;\; {\rm on}  \;\;  G^i\setminus R(a_i^k,b_i^k,l_i),$$
where $C_0>0$ is an absolute constant. By the maximum principle for harmonic functions $p_{\mu_k}(z)>-2C\cos (2\pi ( x-(b_i^k+a_i^k)/2)/l_i)e^{-y/l_i}$, where $z=(x,y)$ and we get the desirable inequality.

 (17) We have proved a lower estimate for the function $p_{\mu_k}$ on the set $\partial G_k $. In the following sections 
  we shall get an upper estimate for the function $u_{I_n}$ on this set and, as a consequence,  the convergence demanded in (15).

(18) First we prove 
$$ \lim \inf v(I_n) > -\infty$$
as $n\to \infty$.

Assume the opposite, i.e.,
$$ \lim \inf v(I_n) = -\infty .$$

Denote $l_i^k=b_i^k-a_i^k$,  $G'_n:=\{ z\in {\mathbb R}^2, p_{\mu_n}(z)<-C_0\}$.
From (16) it follows that 
$$ G'_k \subset \cup_i R(a_i^k, b_i^k, l_i^k) . $$.

From our assumption and (11) it follows

$$ \lim v( \cup_i R(a_i^k,b_i^k,l_i^k)) = -\infty$$
and hence

$$ \lim v(G'_n) = -\infty .$$

 Since $v(G_n)=0$ and $v(G'_n)=v(G_n) +C_0=C_0$,  the result follows.

(19) Let $U(z)$ be the normalized potential of the equilibrium measure of the segment $e=[ -1, 1]$ on the real axis. Then
$$U(z) <C+ \ln |z| $$
near infinity and
$$ U(z) < C\sqrt { {\rm dist} (z,e) } $$
near the segment $e$. The last inequality is well known and it follows for example from the explicit expression  for the 
equilibrium measure of the segment  
$[a, b]$ (cf. \cite{Se}, A2.4)
$$\mu_{[a, b]}=\frac{1}{\pi} \frac{dx}{\sqrt{(b-x)(x-a)}} \; .$$  

We will use the last inequalities for potential $U$ to get upper bounds for the potentials $u_{I_k}(z)$.

Fix $k$ and let $l_i=b_i^k-a_i^k$. Denote by $c_i$ the centre of the segment $[a_i,b_i]$, $c_i=((a_i+b_i)/2,0)$. Denote by $H$ an annulus centered at $c_i$
with inner radius $l_i$ and exterior radius $2R$. Then as it follows from (18),  $u_{I_k}(z)< C$ on the exterior boundary.
Thus
$$u_{I_k}(z) < CU((z-c_i)/l_i) /\ln (1/l_i ) \; .$$

The above inequalities for the function $U$ imply the following estimates.  Fix  $\epsilon >0$. 

(a) For $l_i>\epsilon$ the normalized potential $u_{I_k}(z) <C_\epsilon \sqrt {\epsilon /l_i} / \ln {(1/l_i)}$ on $R(a_i^k,b_i^k,\epsilon) $, where $C>0$ is an absolute constant.  

(b) If $\epsilon / A<l_i\leq \epsilon, \, A>1$ then $u_{I_k}(z) <C_A \ln (A) / \ln ({ 1/l_i})$ on $R(a_i^k,b_i^k,\epsilon) $, where $C_A>0$ is a constant depending on $A$.

(20) In this section we define the sets $H_k$ to verify the inequality of (15). By (13) there is a sequence  $c_k\to 0$ such that
$G_k\subset S_{c_k} $. For any $\epsilon >0$  later on we shall define a sufficiently large constant $C>0$ and a sufficiently small
constant $\delta >0$. Denote $l_i^k=b_i^k-a_i^k$.  For $l_i\leq \delta $ set
$$ R_i^k=  R(a_i^k,b_i^k,Cl_i^k) .$$
 For $l_i> \delta $ set
 $$ R_i^k=  R(a_i^k,b_i^k,d_i^k), $$
 where $d_i^k =\min \{ l_i^k/C, c_k \} $.
Now let
 $$ H_k = \cup_i R_i^k \; .$$
 Then for a sufficiently large constant $C>0$ and a sufficiently small
constant $\delta >0$ by (16) and (19) we get the inequality
 $$| p_{\mu_n}(z)-u_{I_n}(z)|< \epsilon \;\; {\rm on}\; \; \partial H_k \; .$$

(21) We show that   $\Delta p_{\mu_n}(z)- \Delta u_{I_n}(z)\to 0$ as $n\to \infty$
in the sense of Schwartz distributions, hence it converges in the weak topology on Radon measures. Since $\Delta p_{\mu_n}(z)$,  $\Delta u_{I_n}(z)$ are probability measures and since the space of  probability measures on a segment is weakly compact, we can choose a subsequence such that the
following limits exist : 
$\Delta p_{\mu_n}(z)\to \mu_1, \;  \Delta u_{I_n}(z)\to \mu_2$  weakly as $n\to \infty$, where $\mu_1, \mu_2$ are probability measures on $I$. Since $p_{\mu_n}(z)$ is a potential of a probability measure and since weak convergence of measures implies convergence of their potentials, we may  assume  convergence $p_{\mu_n}(z)\to h(z), \, z \in R^2 \setminus I$, where $h$ is a harmonic function in the complement of $I$.  Hence by (20) $u_{I_n}(z)\to h(z)$.  Thus  $\mu_1=\mu_2$. Therefore  $\Delta p_{\mu_n}(z)- \Delta u_{I_n}(z)\to 0$ as $n\to \infty$
in the sense of distributions. 

Denote by $\nu_n$ the equilibrium measure of the set $I_n$. Then $\nu_n= \Delta u_{I_n}$. 
Thus $\nu_n- \Delta p_{\mu_n}(z) \to 0$ weakly.
 Therefore, it follows that  $p_{I_n}(z)\to p_\mu(z)$ weakly, and since $p_{I_n}(z) = u_{I_n}(z) + v(I_n)$ we get $ v(I_n)\to 0$ and
$| p_{\mu_n}(z)-p_{I_n}(z)|\to 0$. We can choose a sequence $\delta_n \to 0$
such that if $I'_n$ is the image of $I_n$  after delation of the axis $x$ to $(1+\delta_n)x$ then $ {v}(I'_n)=0$.
Since  $\delta_n \to 0$, we have $p_{I'_n}(z)\to p_\mu(z)$ weakly. 

The theorem is proved. 

\section {Aknowledgement}
We would like to express our gratitude to J.-P. Serre and S. Vl\v{a}du\c{t} for many valuable remarks.

\end{document}